\documentclass[12pt,a4paper]{article}
\usepackage[utf8x]{inputenc}
\usepackage[english]{babel}
\usepackage[T1]{fontenc}
\usepackage{amsmath}
\usepackage{amsfonts}
\usepackage{amssymb}
\usepackage{makeidx}
\usepackage{graphicx}
\usepackage[left=3cm,right=2cm,top=3cm,bottom=2cm]{geometry}
\usepackage{setspace}
\usepackage{amsthm}
\usepackage{cite}
\usepackage{mathtools, nccmath}

\usepackage[colorlinks=true,urlcolor=blue,
citecolor=red,linkcolor=blue,linktocpage,pdfpagelabels,
bookmarksnumbered,bookmarksopen]{hyperref}

\addto\captionsenglish{}

\newtheorem{definition}{Definition}[section]

\newtheorem{proposition}{Proposition}[section]

\newtheorem{theorem}{Theorem}[section]

\newtheorem{example}{Example}[section]

\newtheorem{lemma}{Lemma}[section]

\newtheorem{observation}{Remark}[section]

\newtheorem{corollary}{Corolary}[section]

\numberwithin{equation}{section}
\newcommand{\cupdot}{\mathbin{\mathaccent\cdot\cup}}

\newcommand{\bo}{\begin{observation}}
\newcommand{\eo}{\end{observation}}
\newcommand{\bd}{\begin{definition}}
\newcommand{\ed}{\end{definition}}
\newcommand{\bp}{\begin{proposition}}
\newcommand{\ep}{\end{proposition}}
\newcommand{\bt}{\begin{theorem}}
\newcommand{\et}{\end{theorem}}
\newcommand{\bc}{\begin{corollary}}
\newcommand{\ec}{\end{corollary}}
\newcommand{\bl}{\begin{lemma}}
\newcommand{\el}{\end{lemma}}
\newcommand{\be}{\begin{example}}
\newcommand{\ee}{\end{example}}
\newcommand{\beq}{\begin{equation}}
\newcommand{\eeq}{\end{equation}}
\newcommand{\beqa}{\begin{equation*}}
\newcommand{\eeqa}{\end{equation*}}
\newcommand{\R}{\mathbb{R}}
\newcommand{\RN}{\mathbb{R}^{N}}
\newcommand{\N}{\mathbb{N}}
\newcommand{\Wsp}{ W^{s,p}(\RN) }
\newcommand{\Lp}{ L^{p}(\mathbb{R}^N) }

\newcommand{\Lpquatro}{L^{\frac{2N}{2N-1}p}(\RN)}

\newcommand{\crit}{\alpha_{s, N}^{\ast}}

\newcommand{\K}{\mathcal{K}}
\newcommand{\A}{\mathcal{A}}

\newcommand{\Ll}{\mathcal{L}}

\newcommand{\Nn}{\mathcal{N}}
\newcommand{\Oo}{\mathcal{O}}

\newcommand{\Wball}{W^{s,p}(B_R)}

\newcommand{\Lomega}{ L^{\omega}(\mathbb{R}^N)}

\newcommand{\intR}{\displaystyle\int\limits_{\mathbb{R}^N}}
\newcommand{\un}{u_{n}}

\newcommand{\uk}{u_{k}}

\newcommand{\vn}{v_{n}}

\newcommand{\yn}{y_{n}}

\newcommand{\until}{\tilde{u}_n}

\newcommand{\RA}{\rightarrow}
\newcommand{\CF}{\rightharpoonup}
\newcommand{\IC}{\hookrightarrow}

\newcommand{\ds}{\displaystyle\int\limits}
\newcommand{\intlog}{ \displaystyle\int\limits_{\mathbb{R}^N} \displaystyle\int\limits_{\mathbb{R}^N} \ln (|x-y|)|u(x)|^p |u(y)|^p dx dy}
\newcommand{\intLog}{ \displaystyle\int\limits_{\mathbb{R}^N} \displaystyle\int\limits_{\mathbb{R}^N} \ln (|x-y|)|u(x)|^p |u(y)|^{p-2}u(y)v(y) dx dy}

\usepackage{times, color, xcolor}
\addto\captionsportuguese{
}

\begin{document}

 	\title{Existence and multiplicity of solutions for the fractional $ p $-Laplacian  Choquard logarithmic 
equation involving a nonlinearity  with  exponential critical and subcritical growth 
		\thanks{The first author was supported    by  Coordination of Superior Level Staff Improvement-(CAPES) -Finance Code 001 and  S\~ao Paulo Research Foundation- (FAPESP), grant $\sharp $ 2019/22531-4,
while the second  author was supported by  National Council for Scientific and Technological Development -(CNPq),   grant $\sharp $ 307061/2018-3 and FAPESP  grant $\sharp $ 2019/24901-3.
			}}
	\author{
		Eduardo  de S. Böer \thanks{Corresponding author} \thanks{ E-mail address: eduardoboer04@gmail.com Tel. +55.51.993673377}  and Ol\'{\i}mpio H. Miyagaki \footnote{ E-mail address: ohmiyagaki@gmail.com, Tel.: +55.16.33519178 (UFSCar).}\\
		{\footnotesize Department of Mathematics, Federal University of S\~ao Carlos,}\\
		{\footnotesize 13565-905 S\~ao Carlos, SP - Brazil}\\ }
\noindent
				
	\maketitle

\noindent \textbf{Abstract:} In the present work we obtain the existence and multiplicity of nontrivial solutions for the Choquard logarithmic equation $(-\Delta)_{p}^{s}u + |u|^{p-2}u + (\ln|\cdot|\ast |u|^{p})|u|^{p-2}u = f(u) \textrm{ \ in \ } \mathbb{R}^N $ , where $ N=sp $, $ s\in (0, 1) $, $ p>2 $, $ a>0 $, $ \lambda >0 $ and $f: \R \RA \R $ a continuous nonlinearity with exponential critical and subcritical growth. We guarantee the existence of a nontrivial solution at the mountain pass level and a nontrivial ground state solution under critical and subcritical growth. Morever, when $ f $ has subcritical growth we prove the existence of infinitely many solutions, via genus theory.

\vspace{0.5 cm}


\section{Introduction}
In the present paper we study existence and multiplicity results for the Choquard logarithmic equation
\begin{equation} \label{P}
(-\Delta)_{p}^{s}u + |u|^{p-2}u + (\ln|\cdot|\ast |u|^{p})|u|^{p-2}u = f(u) \textrm{ \ in \ } \mathbb{R}^N,
\end{equation}
where $ N=sp $, $ s\in (0, 1) $, $ p>2 $, $ a=1 $, $ \lambda =1 $ and $f: \R \RA \R $ is continuous, with primitive $ F(t)=\int\limits_{0}^{t}f(\tau)d\tau $. Precisely, when dealing with critical exponential growth, we will guarantee the existence of a solution at the mountain pass level and a ground state solution, in the sense that it will be the least energy level. On the other hand, considering $ f $ with subcritical exponential growth, we are able to prove that equation (\ref{P}) has infinitely many solutions, via genus theory. Such nonlinearity behaviour is motivated by the Moser-Trudinger Lemma \ref{l2}.

First of all, to make the notation concise, we set, for $ \alpha > 0 $ and $ t \in \R $,
$$ R(\alpha , t) = \exp(\alpha |t|^{\frac{N}{N-s}}) - S_{k_p -2}(\alpha , t) = \sum\limits_{k_p -1}^{+\infty} \dfrac{\alpha^k}{k!} |s|^{\frac{N}{N-s}k},
$$
where $S_{k_p -2}(\alpha , t) = \sum\limits_{k=0}^{k_p -2}\dfrac{\alpha^{k}}{k!}|t|^{\frac{N}{N-s}k}$ and $ k_p = \min\{ k \in \N \ ; \ k \geq p \} $.

Then, we recall that a function $ h $ has \textit{subcritical} exponential growth at $ +\infty $, if
$$
\lim\limits_{t\RA + \infty}\dfrac{h(t)}{R(\alpha , t)} = 0 \textrm{ \ , for all \ } \alpha >0 ,
$$
and we say that $ h $ has \textit{critical} exponential growth at $ +\infty $, if there exists $ \alpha_0 > 0 $ such that
$$
\lim\limits_{t\RA + \infty}\dfrac{h(t)}{R(\alpha , t)} = \left\{ \begin{array}{ll}
0, \ \ \ \forall \ \alpha > \alpha_{0}. \\
+\infty , \ \ \ \forall \ \alpha < \alpha_0 .
\end{array} \right. 
$$

In the following we present some necessary conditions to obtain our main results. This kind of hypothesis are usual in works with Moser-Trudinger inequality, such as \cite{[frac], [9], [boer]}. We assume that $ f $ satisfies
$$ f\in C(\R , \R), f(0)=0, \mbox{has critical exponential growth and} \ F(t) \geq 0 \mbox{ \ , for all \ } t\in \R . \leqno{(f_1)}$$
$$ \lim\limits_{|t|\RA 0} \dfrac{f(t)}{|t|^{p-2}t}=0.  \leqno{(f_2)}$$

From $ (f_1) $ and $(f_2)$, given $ \varepsilon >0 $ and $\alpha > \alpha_0$, fixed, there exists a constant $ b_1> 0$ such that 
\begin{equation}\label{eq1}
|f(t)|\leq \varepsilon |t|^{p-1} +b_1 R(\alpha , t) \ , \ \ \ \forall \ t \in \R .
\end{equation}
As a consequence, 
\begin{equation}\label{eq2}
|F(t)|\leq \dfrac{\varepsilon}{p} |t|^p +b_1|t|R(\alpha , t) \ , \ \ \ \forall \ t \in \R .
\end{equation}

Our strategy to prove Theorem \ref{t1} will consist in finding a Cerami sequence for the mountain pass level. In order to verify that such sequence is bounded in $\Wsp $, we will need the following condition 
$$\mbox{ there exists }\ \theta > 2p \ \mbox{ such that}\   f(t)t \geq \theta F(t) > 0, \ \mbox{for all } \  t> 0 . \leqno{(f_3)}$$

Moreover, as we will be working with an exponential term, to guarantee that the mentioned Cerami sequence and the minimizing sequence for the ground state satisfy the exponential estimates, we rely in the next condition.
$$ \mbox{ there  exist}\  q>2p \ \mbox{ and}\  C_q> \dfrac{[2(q-p)]^{\frac{q-p}{p}}}{q^{\frac{q}{p}}}\dfrac{S_{q}^{q}}{\rho_{0}^{q-p}} \ \mbox{ such that}\  F(t) \geq C_q |t|^q , \ \mbox{for all}\  t\in \R , \leqno{(f_4)}$$ 
for $ S_q, \rho_0 >0 $ to be defined in Lemma \ref{l15}.

Therefore, we are able to state our first main result.

\bt\label{t1}
Assume $ (f_1)-(f_4) $, $ q>2p $ and $ C_q>0 $ sufficiently large. Then,
\begin{itemize}
\item[(i)] Problem (\ref{P}) has a nontrivial solution $ u\in X $ such that $$ I(u)=c_{mp}=\inf\limits_{\gamma \in \Gamma}\max\limits_{t\in [0, 1]}I(\gamma(t)) ,  $$
where $ \Gamma = \{ \gamma \in C([0, 1], X) \ ; \ \gamma(0) \ , \ I(\gamma(1))< 0 \} $. 

\item[(ii)] Problem (\ref{P}) has a nontrivial ground state solution $ u\in X $, that is, $ u $ satisfies
$$ I(u)=c_g = \inf\{ I(v) \ ; \ v\in X \mbox{ \ is a solution of (\ref{P})} \} . $$
\end{itemize}
\et

For the second main result, we are concerned with multiplicity of solutions. However, to obtain this we need to exchange the condition $ (f_1) $ by the condition below.
$$ f\in C(\R , \R), \ f \mbox{ \ is odd, has subcritical exponential growth and} \ F(t) \geq 0 \mbox{, for all \ } t\in \R . \leqno{(f_1 ')}$$
Also, we need to add a condition that gives us the desired geometry for the associated functional. That is, 
$$ \mbox{the function \ } t\mapsto \dfrac{f(t)}{t^{2p-1}} \mbox{ \ is increasing in \ } (0, + \infty) . \leqno{(f_5)}$$
From this condition, since $ f $ is odd, it follows that $ \frac{f(t)}{t^{2p-1}} $ is decreasing in $ (-\infty , 0) $.

Moreover, since we can control the exponent using $ \alpha >0 $, we can change condition $ (f_4) $ by a more general condition, that is
$$\mbox{ there  exists } q>2p \mbox{ \ and \ } M_1 > 0 \ \mbox{ such that}\  F(t) \geq M_1 |t|^q \ , \ \forall \ t\in \R , \leqno{(f_4 ')}$$
 
\bt\label{t2}
Suppose $ (f_1 '), (f_2), (f_3), (f_4 '), (f_5) $. Then, problem (\ref{P}) has infinitely many solutions. 
\et

To conclude this introduction, we make a quick overview of the state of art.  In the first half of this overview, we will be concerned with fractional $ p $-Laplacian equations and problems that deal with exponential nonlinearities.  We recall that problems with nonlocal operators arise in many areas, such as optimization, finance, phase transitions, stratified materials, anomalous diffusion, crystal dislocation, soft thin films, semipermeable membranes, flame propagation, conservation laws and water waves. See e.g. \cite{[hitch]}, where the authors provide an extensive list with references for the mentioned applications. For fractional problems of the form
\begin{equation}\label{eq22}
(-\Delta)_{p}^{s}u + V(x)|u|^{p-2}u = f(x, u) + \varepsilon h(x) \textrm{ \ in \ } \mathbb{R}^N,
\end{equation}
we mention \cite{[w1], [w2], [w4]}. In the \cite{[w1]}, the authors consider problem (\ref{eq22}) in the case $ N=sp $, $ s\in (0, 1) $, $ V(x) $ having a positive lower bound and being coercive or satisfying $ \frac{1}{V(x)}\in L^{1}(\RN) $, $ f(x, t) $ behaving like $ e^{\alpha|t|^{\frac{N}{N-s}}} $ at infinity, $ h\in (\Wsp)^{\ast} $ and $ \varepsilon >0 $. Under suitable conditions over $ V, f, h $, the authors guarantee the existence of weak solutions for (\ref{eq22}). In \cite{[w2]}, they consider  $ V\in C(\RN) $ having a positive lower bound, $ p \geq 2 $, $ s\in (0, 1) $, $ N \geq 2 $, $ h \equiv 0 $ and $ f $ $ p$-superlinear. Using Mountain Pass Theorem for Cerami condition, the authors prove the existence of a nontrivial radially symmetric solution. Finally, in \cite{[w4]}, the authors consider $ V\in L^{\infty}(\RN) $ possibly indefinite, a concave and convex nonlinearity $ f(x, u)= w_1(x)|u|^{q-2}u - w_2(x)|u|^{r-2}u $, for $ 0<s<1<q<p<r $, and $ h \equiv 0 $. They begin showing that $(-\Delta)_{p}^{s}u = \lambda V(x)|u|^{p-2}u $ posses a infinite sequence of eigenvalues and that the first one is simple. Then, using this fact, the obtain the existence of infinitely many solutions for (\ref{eq22}). See also \cite{[w3]}.

We also refer the reader for the works \cite{[9], [weyl], [pucci1], [pucci2]}, in which the authors dealt with fractional Laplacian operator, and  the works \cite{[5], [Lam], [17]}, for general problems with Moser-Trudinger type behaviour.

In the second half, we take a look into works that deal with Choquard logarithmic equations. We can cite the recent works of \cite{[6], [cjj], [10], [alves], [wen], [guo]}, where the authors study the following class of equation
\beq\label{i2}
-\Delta u + V(x) u + \gamma (\ln |\cdot| \ast |u|^2) u = f(u), \textrm{ \  \ in \ } \mathbb{R}^N ,
\eeq
under distinct conditions. We briefly discuss some of these works and left others to the reader. In \cite{[6]}, the authors have proved the existence of infinitely many geometrically distinct solutions and a ground state solution, considering $ N=2 $, $ V: \RN \RA (0, \infty) $ continuous and $ \mathbb{Z}^2 $-periodic, $ \gamma > 0 $, $ b\geq 0 $, $ f(u) = b|u|^{p-2}u $ and $ p\geq 4 $. Here  because of the  periodic setting, the global Palais–Smale condition can fail, since  the corresponding functional become invariant under $ \mathbb{Z}^2 $-translations. Then, intending to fill the gap, in \cite{[10]} the authors proved the existence of a mountain pass solution and a ground state solution for the equation (\ref{i2}) in the case $ N=2 $, $ V(x)\equiv a > 0 $, $ \lambda >0 $, $ f(u) = b|u|^{p-2}u $ and $ 2<p<4 $. Also, they verified that, if $ p\geq 3 $, both levels are equal and provided a characterization for them.  In \cite{[cjj]}, the authors dealt with the existence of stationary waves with prescribed norm considering $ \lambda \in \R $. Finally, in \cite{[alves]}, the authors proved the existence of a ground state solution for equation (\ref{i2}), with a nonlinearity of Moser-Trudinger type. We also refer to \cite{[4], [14]} for Choquard equations.

We also call attention to \cite{[boer]}, where the authors consider the following equation
\begin{equation} \label{eq23}
(-\Delta)^{\frac{1}{2}} u + u + (\ln|\cdot|\ast |u|^{2})u = f(u) \textrm{ \ in \ } \mathbb{R},
\end{equation}
with $ f\in C(\R , \R) $ a nonlinearity of Moser-Trudinger type. In this work, they prove that (\ref{eq23}) has a mountain pass solution and a ground state solution, when $ f $ has exponential critical behaviour. On the other, consider $ f $ odd and with subcritical exponential growth, the authors obtained infinitely many solutions for (\ref{eq23}). How to obtain multiplicity results for (\ref{eq23}) in the critical case remains an open problem.

The present work aims to expand and complement those results already found in the literature, considering a Choquard logarithmic equation for the fractional $ p $-Laplacian operator with a nonlinearity under exponential growth. 

Throughout the paper, we will use the following notations: $ \Lomega $ denotes the usual Lebesgue space with norm $ ||\cdot ||_\omega $ \ ; \ $ X' $ denotes the dual space of $ X $ \ ; \ $ B_r(x) $ is the ball in $ \R $ centred in $ x $ with radius $ r>0 $ and simply $ B_r $ when $ x=0 $ \ ;  \ $ C, C_1, C_2, ... $ will denote different positive constants whose exact values are not essential to the exposition of arguments. 

The paper is organized as follows: in section 2 we present the framework's problem and some technical and essential results, mostly concerning the associated functional. Section 3 consists in the proof of a key proposition and our first main result. Finally, in section 4, we prove our second main result, in which we deal with the subcritical case and obtain multiplicity. 

\section{Framework for Problem (\ref{P})}

In this section we establish the necessary framework for solving (\ref{P}). For each $ s \in (0, 1) $ and $ p>2 $, we consider the Sobolev space $ \Wsp = \{ u \in \Lp \ ; \ [u]_{s, p} < + \infty \} $, where
$$
[u]_{s, p}^{p}= \intR \intR \dfrac{|u(x)-u(y)|^p}{|x-y|^{N+ps}} dx dy ,
$$
is the Gagliardo seminorm. It is well-known that the space $ (\Wsp , ||\cdot||) $, where $ ||\cdot||^{p}= [\cdot]_{s, p}^{p} + ||\cdot||_{p}^{p} $, is an uniformly convex Banach space, particularly reflexive, and separable (see \cite[Theorem A.3]{[w5]}). 

Next, we recall that, for a function $ u\in C_{0}^{\infty}(\RN) $, the fractional $ p $-Laplacian operator $ (-\Delta)_{p}^{s} $ is given by
$$
(-\Delta)_{p}^{s}u(x) = C(N, s) \lim\limits_{\varepsilon \RA 0} \ds_{\RN \setminus B_\varepsilon (x)}\dfrac{|u(x)-u(y)|^{p-2}(u(x)-u(y))}{|x-y|^{N+ps}} dy \ , \ \forall \ x \in \RN ,
$$
where the normalizing constant $ C(N, s) $ is defined in \cite{[hitch]}. Throughout this paper we omit the normalizing constant to simplify the expressions. We also remember the reader that $ C_{0}^{\infty}(\RN) $ is dense in $ \Wsp $ (see \cite[Theorem 7.38]{[adams]}). For further considerations about $ \Wsp $ and $ (-\Delta)_{p}^{s} $ and some useful results, we refer to \cite{[hitch], [demengel], [adams], [weyl], [pucci1], [pucci2]}.

We recall that $ u\in \Wsp $ is a weak solution for (\ref{P}) if
\begin{align}
& \intR \intR \dfrac{|u(x)-u(y)|^{p-2}(u(x)-u(y))(v(x)-v(y))}{|x-y|^{N+sp}} dx dy + a \intR |u|^{p-2}uv dx \nonumber \\
& + \intLog = \intR f(u) v dx \ , \ \forall \ v \in \Wsp . \label{eq3}
\end{align}
In this sense, we consider the associated functional $ I : \Wsp \RA \R \cup \{+ \infty\} $ given by
\begin{equation}\label{eq4}
I(u) = \dfrac{1}{p}||u||^p + \dfrac{1}{2p}\intlog - \intR F(u) dx .
\end{equation}
Then, as we will see, critical points of $ I $ will be weak solutions for (\ref{P}). However, one can see that $ I $ is not well-defined over the whole space $ \Wsp $, since we have the term involving the logarithmic function. Hence, following the ideas introduced by Stubbe \cite{[21]}, we consider the slightly smaller space
$$
X = \left\{ u \in \Wsp \ ; \ \intR \ln(1+|x|)|u(x)|^p dx < + \infty \right\} .
$$ 
It is possible to verify that
$$
||u||_\ast = \intR \ln(1+|x|)|u(x)|^p dx 
$$
defines a norm in $ X $. The next natural step is to guarantee that $ X $ enjoys of all good properties that we need. So, we start defining a measure $ \eta : \Ll \RA [-\infty , + \infty] $ given by 
$$
\eta(E)=\ds_{E}\ln(1+|x|)|u(x)|^p dx ,
$$
where $ \Ll $ is the Lebesgue $ \sigma $-algebra in $ \RN $. Then, we can consider the measure space $ L^{p}(\RN , \Ll , \eta) $ and, for each $ u: \RN \RA \R $ Lebesgue measurable, we have $ ||u||_{L^{p}(\RN , \Ll , \eta)} = ||u||_\ast $. 

Consequently, from Riesz-Fischer Theorem (see Folland \cite{[folland]}), the definition of $ L^p $ spaces and the fact that $ X = \Wsp \cap L^{p}(\RN , \Ll , \eta) $, we conclude that $ X $ is a Banach space endowed with the norm $ ||\cdot||_{X}^{p}=||\cdot||^p + ||\cdot||_{\ast}^{p} $. Moreover, as \cite[Proposition A.6]{[pucci1]}, one can show that $ X $ is uniformly convex. Hence, from \cite[Theorem 3.31]{[brezis]}, $ X $ is reflexive. 

The following lemma plays a key role on the continuity of the functional $I$ and allow us to verify that $I$ is lower semicontinuous for $ \Wsp $ and weakly lower semicontinuous for $ X $. In order to do that, we need a result in order to control the exponential term. 

\bl\label{l1}
Let $ (\varphi_n) \subset X $ and $ \varphi \in X $. Then, 
\begin{itemize}
\item[(a)] if $ \varphi_n \RA \varphi $ in $ \Wsp $ or $ \varphi_n \CF \varphi $ in $ X $, then there exists a subsequence $ (\varphi_{n_k})\subset (\varphi_n) $ and a function $ h\in \Wsp $ such that  $ \varphi_{n_k}(x) \RA \varphi(x) $ a.e. in $ \RN $ and $ |\varphi_{n_k}(x)|\leq h(x) $, for all $ k\in \N $ and a.e. in $ \RN $.

\item[(b)] if $ \varphi_n \RA \varphi $ in $ X $, then there exists a subsequence $ (\varphi_{n_k})\subset (\varphi_n) $ and a function $ h\in X $ such that  $ \varphi_{n_k}(x) \RA \varphi(x) $ a.e. in $ \RN $ and $ |\varphi_{n_k}(x)|\leq h(x) $, for all $ k\in \N $ and a.e. in $ \RN $.
\end{itemize}
\el 
\begin{proof}
We will use the construction done in \cite[Proposition 2.7]{[n1]}. For items (a) and (b), similarly as \cite[Lemma 2.4]{[boer]}, one can see that the function $ w $ constructed in \cite{[n1]} belongs to $ \Wsp $. Here, we only highlight, why in item (b) we have $ w \in X $.

By construction, one has $ ||w_n||_{\ast} \leq 1 $, for all $ n \in \N $. Since $ w_{n+1}(x) \geq w_n (x) $ and $ w_n(x)\RA w(x) $ a.e. in $ \RN $, follows $ \ln(1+|x|)|g_j(x)|^p \RA \ln(1+|x|)|g(x)|^p  $ a.e. in $ \RN $ and, from the Monotone Convergence Theorem, 
$$
\intR \ln(1+|x|)|w(x)|^p dx = \lim \intR \ln(1+|x|)|w_n(x)|^p dx \leq 1 < + \infty .
$$
Therefore, $ w\in X $. 
\end{proof}

We would like to call attention that, essentially, the above lemma tell us that the obtained function $ h $ satisfies $ ||h||_\ast < + \infty $. It will be very important in order to verify that $ I $ is $ C^1 $. 

To finish this part, we consider the operators $ A: \Wsp \RA (\Wsp)^{\ast} $ given by
$$
A(u)(v) = \intR \intR \dfrac{|u(x)-u(y)|^{p-2}(u(x)-u(y))(v(x)-v(y))}{|x-y|^{N+sp}} dx dy \ , \ \forall \ u, v \in X
$$
and $ \tilde{A}: \Wsp \RA (\Wsp)^{\ast} $ defined as
$$
\tilde{A}(u)(v) = A(u)(v) + \intR a|u|^{p-2}u v dx \ , \ \forall \ u, v \in X .
$$
One can easily verify that $ \tilde{A}(u)(v) \leq ||u||^{p-1}||v|| $ and $ \tilde{A}(u)(u) = ||u||^p $, for all $ u, v \in \Wsp $. Therefore, as $ \Wsp $ is uniformly convex, from \cite[Proposition 1.3]{[ravi]}, $ \tilde{A} $ verifies the $ (S) $ property, that is, for any sequence $ (\un)\subset \Wsp $ satisfying $ \un \CF u $ in $ \Wsp $ and $ \tilde{A}(\un)(\un-u)\RA 0 $, there exists a subsequence, still denoted by $ (\un) $, such that $ \un \RA u $ in $ \Wsp $.

For the second half of this section, we present some embedding results for $ X $ and verify that $ I $ is well-defined over $ X $ and of class $ C^1 $. Clearly, $ X \IC \Wsp $, once $ ||\cdot||\leq ||\cdot||_X $. Moreover, as the proof of the compactness embedding can be done similarly as \cite[Lemma 2.1]{[boer]}, we omit it here. 

\bp\label{p1}
The space $ X $ is compactly embedded in $ \Lomega $, for all $ \omega \geq p $.
\ep

Before our next result, we recall the celebrated Moser-Trudinger Lemma for unbounded domains and a very useful lemma concerning $ R(\alpha , t) $.

\bl\label{l2}
(Moser-Trudinger Lemma \cite[Theorem 1.1]{[w1]}) Let $ s\in (0, 1) $ and $ sp=N $. Then, there exists $ \crit >0 $ such that for every $ 0\leq  \alpha < \crit $, the following inequality holds
$$
\sup\limits_{u\in \Wsp , ||u||\leq 1} \intR R(\alpha, u) dx < + \infty .
$$
\el

\bl\label{l3}
(\cite[Lemma 2.3]{[w6]}) Let $ \alpha > 0 $ and $ r>1 $. Then, for every $ \beta > r $, there exists a constant $ C_\beta = C(\beta) > 0 $ such that
$$
(\exp(\alpha |t|^{p'}) - S_{k_p -2}(\alpha ,t))^r \leq C_\beta (\exp(\beta \alpha |t|^{p'} - S_{k_p -2}(\beta \alpha , t)) ,
$$
with $ \frac{1}{p}+\frac{1}{p'}=1 $. 
\el

\bl\label{l4}
Let $ \alpha > 0 $. Then, $ R(\alpha , u) \in L^{1}(\RN) $, for all $ u\in \Wsp $.
\el
\begin{proof}
Let $ u\in \Wsp \setminus \{0\} $ and $ \varepsilon > 0 $. Since $ C_{0}^{\infty}(\RN) $ is dense in $ \Wsp $, there exists $ \phi \in C_{0}^{\infty}(\RN) $ such that $ ||u - \phi||< \varepsilon $. Observe that, for each $ k \geq k_p - 1 $,
$$
|u|^{\frac{N}{N-s}k} \leq 2^{\frac{N}{N-s}k}\varepsilon^{\frac{N}{N-s}k} \left|\dfrac{u-\phi}{||u-\phi||}\right|^{\frac{N}{N-s}k}+ 2^{\frac{N}{N-s}k}|\phi|^{\frac{N}{N-s}k} .
$$
Consequently,
$$
R(\alpha, u) \leq R\left( \alpha 2^{\frac{N}{N-s}}\varepsilon^{\frac{N}{N-s}} , \left|\dfrac{u-\phi}{||u-\phi||}\right|^{\frac{N}{N-s}}\right) + R(\alpha 2^{\frac{N}{N-s}} , |\phi|^{\frac{N}{N-s}}) .
$$
From Lemma \ref{l2}, choosing $ \varepsilon > 0 $ sufficiently small such that $ \alpha 2^{\frac{N}{N-s}}\varepsilon^{\frac{N}{N-s}} < \crit $, we have
$$
\intR R\left( \alpha 2^{\frac{N}{N-s}}\varepsilon^{\frac{N}{N-s}} , \left|\dfrac{u-\phi}{||u-\phi||}\right|^{\frac{N}{N-s}}\right) dx < + \infty .
$$
On the other side, since $ \exp(\alpha 2^{\frac{N}{N-s}}|\phi|^{\frac{N}{N-s}}) = \sum\limits_{k=0}^{+ \infty} \dfrac{\alpha^k}{k!}2^{\frac{N}{N-s}k}|\phi|^{\frac{N}{N-s}k} $, there exists $ k_0 \in \N $ such that $ \sum\limits_{k=k_0}^{+ \infty} \dfrac{\alpha^k}{k!}2^{\frac{N}{N-s}k}|\phi|^{\frac{N}{N-s}k} < \varepsilon $. This fact, combined with the fact that $ \frac{N}{N-s}k > 0 $ for all $ k_p -1 \leq k \leq k_0 $, give us
$$
\intR R(\alpha 2^{\frac{N}{N-s}} , |\phi|^{\frac{N}{N-s}}) dx = \ds_{supp \phi} R(\alpha 2^{\frac{N}{N-s}} , |\phi|^{\frac{N}{N-s}}) dx < + \infty .
$$
Therefore, $ R(\alpha, u) \in L^1 (\RN) $, for all $ u\in \Wsp $.
\end{proof}

\bo\label{obs1}
From Lemmas \ref{l3} and \ref{l4}, we conclude that $ R(\alpha , u)^l \in L^{1}(\RN) $, for all $ u\in \Wsp $, $ \alpha > 0 $ and $ l \geq 1 $.
\eo

\bo\label{obs2}
Given $ \varepsilon > 0 $, from $ (f_3) $, there exists $ \delta > 0 $ such that $ |f(u)|\leq \varepsilon |u|^{p-1} $, for all $ |u|\leq \delta $. Now, as $ q> p $, there exists $ r> 0 $ such that $ q = p + r $. So, once $ z^{q-1}, z^{\frac{N}{N-s}k} $, for all $ k \geq k_p -1 $, and $ z^r $ are increasing functions, for $ |u|\geq \delta $,
$$
|F(u)| \leq \dfrac{\varepsilon}{p}|u|^p + b_1 |u|R(\alpha , u) \leq \dfrac{\varepsilon}{p}|u|^p \dfrac{|u|^r R(\alpha , u)}{\delta^r R(\alpha, \delta)} + b_1 |u| \dfrac{|u|^{q-1}}{\delta^{q-1}} R(\alpha , u) = b_2 |u|^q R(\alpha , u) ,
$$
where $ b_2 = \frac{\varepsilon}{p \delta^r R(\alpha, \delta)}+ \frac{b_1}{\delta^{q-1}}>0 $. Therefore, for $ \alpha > \alpha_0 $,
\begin{equation}\label{eq5}
|F(u)| \leq \dfrac{\varepsilon}{p}|u|^p + b_2 |u|^q R(\alpha , u)  \ , \ \forall \ u \in X .
\end{equation}
\eo

Let $ t, t' > 1 $, with $ \frac{1}{t}+\frac{1}{t'}=1 $. From (\ref{eq5}), Remark \ref{obs1}, $ X \IC \Wsp \IC \Lomega $, for all $ \omega \geq p $, and Hölder inequality, we conclude that
\begin{equation}\label{eq6}
\intR |F(u)| dx \leq \dfrac{\varepsilon}{p}||u||_{p}^{p}+b_2 ||u||_{qt}^{q}\left( \intR R(\alpha, u)^{t'} dx \right)^{\frac{1}{t'}} < + \infty \ , \ \forall \ u \in X .
\end{equation}

\bo
Analogously as in Remark \ref{obs2}, for every $ \varepsilon > 0 $, $ \alpha > \alpha_0 $, $ q > p $ and $ u \in X $, one can obtain
\begin{equation}\label{eq7}
|f(u)| \leq \varepsilon |u|^{p-1}+b_1 |u|^{q-1}R(\alpha , u) .
\end{equation}
\eo

Next, inspired by \cite{[6]}, we define three auxiliar functionals $V_1: \Wsp \RA [0, \infty],$ $V_2: L^{\frac{2N}{2N-1}p}(\RN) \RA [0, \infty) $ and $V_0: \Wsp \RA \R \cup \{\infty\},$ given by 
\beqa
u \mapsto V_1(u)=\intR \intR \ln(1+|x-y|)|u(x)|^p |u(y)|^p dx dy ,
\eeqa
\beqa
u \mapsto V_2(u, v)=\intR \intR \ln\left(1+\dfrac{1}{|x-y|}\right)|u(x)|^p |u(y)|^p dx dy ,
\eeqa
\beqa
u \mapsto V_0(u, v)=V_1(u, v)-V_2(u, v)=\intR \intR \ln(|x-y|)|u(x)|^p |u(y)|^p dx dy. 
\eeqa
These definitions are understood to being over measurable function $u, v: \RN \RA \R $, such that the integrals are defined in the Lebesgue sense. 

\bo\label{obs4}
\noindent \textbf{(i)} As a consequence of Hardy-Littlewood-Sobolev Inequality (HLS) \cite{[15]}, with $ \alpha = \beta = 0 $ and $ \lambda = 1 $, we have $ \frac{1}{q}+\frac{1}{t}+\frac{1}{N}=2 $. So, making a natural choice for $ q $ and $ t $, that is $ q=t=\frac{2N}{2N-1} $, we obtain
\beq \label{v2}
|V_2(u)|\leq K_0||u||_{\frac{2N}{2N-1}p}^{2p} \ , \ \ \forall \ u\in L^{\frac{2N}{2N-1}p}(\RN) ,
\eeq
so $ V_2 $ takes finite values over $ L^{\frac{2N}{2N-1}p}(\RN)$. 

\noindent \textbf{(ii)} $ V_1(u)\leq 2||u||_{\ast}^{p}||u||_{p}^{p} $, since $ \ln(1+|x-y|)\leq \ln(1+|x|)+\ln(1+|y|) $.

\noindent \textbf{(iii)} $\intR \intR \ln(1+|x-y|)|u(x)|^p |v(y)|^p dxdy \leq ||u||_{\ast}^{p}||v||_{p}^{p}+ ||v||_{\ast}^{p}||u||_{p}^{p}$.
\eo

\bl\label{l25}
Let $ (\un) \subset X $ and $ u\in X $ such that $ \un \RA u $ on $ \Wsp $. Then, we have
$$
\intR F(\un) \RA \intR F(u) \ \ , \ \ \intR f(\un)\un \RA \intR f(u)u \mbox{ \ \ and \ \ } \intR f(\un) v \RA \intR f(u)v \ , \ \forall \ v \in X .
$$
\el
\begin{proof}
Since $ \un \RA u $  in $ \Wsp $, $ \un(x) \RA u(x) $ a.e. in $ \R $ and, from \cite[Theorem 6.9]{[hitch]}, $ \un \RA u $ in $ \Lomega $ for all $ \omega \geq p $. By Lemma \ref{l1} and using the Dominated Convergence Theorem, the result follows.
\end{proof}

One can see, from (\ref{eq6}) and Remark \ref{obs4}, that $ I(u)< + \infty $ for all $ u\in X $.

\bl\label{c1}
The functionals $ V_1 , V_2 , V_0 , I$ are of class $ C^1 (X, \R) $, with
\begin{equation}\label{eq8}
V_{1}'(u)(v)=2p\intR \intR \ln(1+|x-y|)|u(x)|^p |u(y)|^{p-2}u(y)v(y) dxdy .
\end{equation}
and
\begin{equation}\label{eq9}
V_{2}'(u)(v)=2p\intR \intR \ln(1+\dfrac{1}{|x-y|})|u(x)|^p |u(y)|^{p-2}u(y)v(y) dxdy .
\end{equation}
\el
\begin{proof}
The proof follows by a standard way, combining the Mean Value Theorem with Lemma \ref{l1}, for $ V_1 $ and $ V_2 $, and  by Lemma \ref{l25} and \cite[Lemma 2]{[pucci2]}.
\end{proof}

\bo\label{obs5}
Note that, if $ u\in X $, then $ u\in \Lpquatro $. So, one can see that $ V_2 $ is of class $ C^1 $ in $ \Lpquatro $.
\eo

\bl\label{l24}
\begin{itemize}
\item[(i)] The functional $ V_1 $ is weakly lower semicontinuous in $ \Wsp $. 

\item[(ii)] The functional $ I $ is weakly lower semicontinuous in $ X $.

\item[(iii)] The functional $ I $ is lower semicontinuous in $ \Wsp $. 
\end{itemize}
\el
\begin{proof}
\textbf{(i)} Let $ (\un)\subset \Wsp $ such that $ \un \CF u $. Then, $ \un \CF u $ in $ \Wball $ and, from \cite[Theorem 7.1]{[hitch]}, $ \un \RA u $ in $ L^p(B_R) $, for all $ R>0 $. We claim that
$$ \lim\limits_{n\RA + \infty} \ds_{B_R} \ds_{B_R} \ln(1+|x-y|)|\un(x)|^p|\un(y)|^p dx dy = \ds_{B_R} \ds_{B_R} \ln(1+|x-y|)|u(x)|^p|u(y)|^p dx dy .$$
Indeed, up to subsequence, one can see that 
$
\ds_{B_R}||\un(y)|^p - |u(y)|^p| dy \RA 0 .
$
Then, 
\begin{align*}
& \left| \ds_{B_R} \ds_{B_R} \ln(1+|x-y|)|\un(x)|^p|\un(y)|^p dx dy - \ds_{B_R} \ds_{B_R} \ln(1+|x-y|)|u(x)|^p|u(y)|^p dx dy \right| \\
& \leq [\ln(1+2R)||\un||_{p}^{p} + \ln(1+2R)||u||_{p}^{p}]\ds_{B_R}||\un(x)|^p - |u(x)|^p| dx \RA 0 .
\end{align*}
As a consequence, for each $ R>0 $, 
$$
\liminf V_1(\un) \geq \ds_{B_R} \ds_{B_R} \ln(1+|x-y|)|u(x)|^p|u(y)|^p dx dy .
$$
Hence, from the Monotone Convergence Theorem,
$$
\liminf V_1(\un) \geq \lim\limits_{R \RA + \infty} \ds_{B_R} \ds_{B_R} \ln(1+|x-y|)|u(x)|^p|u(y)|^p dx dy =V_1(u) .
$$

\noindent \textbf{(ii)} Follows from the fact that $ V_2 \in C^{1}(L^{\frac{2N}{2N-1}p}, \R) $, (\ref{v2}), $ X \hookrightarrow \Wsp $, Proposition \ref{p1} and item (i).

\noindent \textbf{(iii)} Follows from (\ref{v2}), \cite[Theorem 6.9]{[hitch]}, item (i) and the fact that $ I(u)-\frac{1}{2p}V_1(u) $ is continuous with respect to $ ||\cdot|| $.
\end{proof}

\section{Geometry of $ I $ and Technical Results}

In this section we will provide some technical results and investigate the geometry of $I$. First of
all, we will verify some conditions that can allow us to get convergence in $X$. Then, we see how
one can get boundedness for the exponential term. Finally, we verify that $I$ has the mountain
pass geometry.

\bl\label{l9}
Let $ (\un)\subset X $ such that $ \un \rightharpoonup u $ in $ X $. Then, 
$$
\lim\limits_{n\RA + \infty} \intR \intR \ln(1+|x-y|)|\un(x)|^p |u(y)|^{p-2} u(y) (\un(y) - u(y)) dx dy = 0 .
$$
\el
\begin{proof}
First of all, we have
\begin{align*}
& \left|\intR \intR \ln(1+|x-y|)|\un(x)|^p |u(y)|^{p-2} u(y) (\un(y) - u(y)) dx dy\right| \\
& \leq ||\un||_{\ast}^{p}||u||_{p}^{p-1}||\un - u||_{p} + \intR \intR \ln(1+|y|)|\un(x)|^p |u(y)|^{p-1}|\un(y) - u(y)| dx dy .
\end{align*}
Now, for any $ R>0 $, fixed, define
$$
\intR \ln(1+|y|)|u(y)|^{p-1}|\un(y) - u(y)| dy = h_n(R)+g_n(R) ,
$$
where
$$
h_n(R) = \int\limits_{B_R} \ln(1+|y|)|u(y)|^{p-1}|\un(y) - u(y)| dy
$$
and 
$$
g_n(R) = \int\limits_{B_{R}^{c}} \ln(1+|y|)|u(y)|^{p-1}|\un(y) - u(y)| dy .
$$
For $ y\in B_R $, $ \ln(1+|y|) \leq \ln(1+R) $. Since $ \un \rightharpoonup u $ in $ X $, by Proposition \ref{p1}, $ \un \RA u $ in $ \Lomega $, for all $ \omega \geq p $. Thus,
$$
|h_n(R)| \leq \ln(1+R)||u||_{p}^{p-1}||\un - u||_p \RA 0 ,
$$
as $ n\RA + \infty $. On the other hand, 
\begin{equation}\label{eq10}
g_n(R) \leq (\int\limits_{B_{R}^{c}} \ln(1+|y|)|u(y)|^{p}dy)^{\frac{1}{p'}}(\int\limits_{B_{R}^{c}} \ln(1+|y|)|\un(y) - u(y)|^p dy)^{\frac{1}{p}} .
\end{equation}
Note that 
$$
(\int\limits_{B_{R}^{c}} \ln(1+|y|)|\un(y) - u(y)|^p dy)^{\frac{1}{p}} \leq C_1 (||\un||_\ast + ||u||_\ast ) \leq C_2 ,
$$
since $ (\un)\subset X $ is bounded. Then, from (\ref{eq10}),
$$
g_n(R) \leq C_1 \left(\int\limits_{B_{R}^{c}} \ln(1+|y|)|u(y)|^{p}dy\right)^{\frac{1}{p'}} = C_1 \varphi(R) \RA 0 \mbox{ \ ,  as \ } R \RA + \infty .
$$
Consequently, for all $ R> 0 $, 
$$
\limsup\limits_{n\RA + \infty} \left| \intR \intR \ln(1+|x-y|)|\un(x)|^p |u(y)|^{p-2} u(y) (\un(y) - u(y)) dx dy \right| \leq C_1 \varphi(R) .
$$
Thus, taking $ R\RA + \infty $, $ \varphi(R) \RA 0 $ and we conclude the proof.
\end{proof}

In order to prove the next proposition, we will need the following technical lemma, which is, essentially, a corollary of Ergorov's Theorem.

\bl\label{l10} 
Let $ u\in \Lp \setminus \{0\} $ and $ (\un)\subset \Lp $ such that $ \un(x)\RA u(x) $ a.e in $ \RN $. Then, there exists $ R\in \N $, $ \delta > 0 $, $ n_0 \in \N $ and $ A\subset B_R $, such that $ A $ is measurable, $ \mu(A)>0 $ and $ \un(x) > \delta $, for all $ x\in A $ and for all $ n \geq n_0 $.
\el

\bp\label{p2}
Let $ u\in \Lp \setminus \{0\} $, $ (\un)\subset \Lp $ such that $ \un(x)\RA u(x) $ a.e in $ \RN $ and $ (\vn)\subset \Lp $ bounded. If
\begin{equation}\label{eq11}
\alpha = \sup\limits_{n} \intR \intR \ln(1+|x-y|)|\un(x)|^p |\vn(y)|^p dx dy < +\infty ,
\end{equation}
then $ ||\vn||_\ast $ is bounded. Moreover, setting
$$
\alpha_n = \intR \intR \ln(1+|x-y|)|\un(x)|^p |\vn(y)|^p dx dy ,
$$
for each $ n\in \N $, if $ \alpha_n \RA 0 $ and $ ||\vn||_p \RA 0 $, then $ ||\vn||_\ast \RA 0 $.
\ep
\begin{proof}
Let $ n_0, R, \delta $ and $ A $ is in Lemma \ref{l10}. Then, $ \un(x)> \delta $, for all $ n\geq n_0 $. From $ \ln $ properties and the fact that $ \alpha_n \geq 0 $, for all $ n\in \N $, we have
\begin{align*}
\alpha_n & \geq \int\limits_{B_{2R}^{c}}\int\limits_{A} \ln(1+|x-y|)|\un(x)|^p |v(y)|^p dxdy \\
& > \delta^p \int\limits_{B_{2R}^{c}}\int\limits_{A} \ln(\sqrt{1+|y|})|\vn(y)|^p dx dy \\
& = \dfrac{\delta^p \mu(A)}{2} \left( \intR\ln(1+|y|)|\vn(y)|^p dy - \int\limits_{B_{2R}}\ln(1+|y|)|\vn(y)|^p dy \right) \\
& \geq \dfrac{\delta^p \mu(A)}{2}  (||v||_{\ast}^{p}-\ln(1+2R)||\vn||_{p}^{p}).
\end{align*}
Consequently, 
\begin{equation}\label{eq13}
0\leq ||\vn||_{\ast} \leq \left(\dfrac{2}{\mu(A)\delta^p}\alpha_n + \ln(1+2R)||\vn||_{p}^{p}\right)^{\frac{1}{p}} .
\end{equation}
Therefore, from hypothesis, $ (||\vn||_\ast)\subset \R $ is bounded. Moreover, if $ \alpha_n \RA 0 $ and $ ||\vn||_p \RA 0 $, from equation (\ref{eq13}), $ ||\vn||_\ast \RA 0 $. 
\end{proof}

\bl\label{l11}
Let $ u\in \Wsp $, $ r>p $, $ l\geq 1 $, $ \beta > 0 $ and $ ||u||\leq M $, for $ M>0 $ sufficiently small. Then, there exists a constant $ K_1 = K_1(\beta , N, M, l, s) > 0 $ such that
$$
\intR |u|^r R(\beta , u)^l dx \leq K_1 ||u||_{t_0}^{r} ,
$$
for some $ t_0 > p $.
\el
\begin{proof}
From Lemma \ref{l3}, for $ \beta_1 = \beta_1(l) > l $ with $ \beta_1 \sim l $, there exists a constant $ C_1 =C_1(\beta_1) > 0 $ such that $ R(\beta , u)^l \leq C_1 R(\beta_1 \beta , u) $. 

Let $ t, t' >1 $ with $ \frac{1}{t}+\frac{1}{t'}=1 $. From Hölder inequality, 
$$
\intR |u|^r R(\beta , u)^l dx \leq C_1 \intR |u|^r R(\beta_1 \beta , u) dx \leq C_1 \left(\intR  R(\beta_1 \beta , u)^t dx\right)^{\frac{1}{t}}||u||_{rt'}^{r} .
$$
Once again, from Lemma \ref{l3}, for $ \beta_2 = \beta_2(t) > t $ with $ \beta_2 \sim t $, there exists a constant $ C_2=C_2(\beta_2) >0 $ satisfying $ R(\beta_1 \beta , u)^t \leq C_2 R(\beta_2 \beta_2 \beta , u) $. Thus,
$$
\intR |u|^r R(\beta , u)^l dx \leq C_1 C_{2}^{\frac{1}{t}} \left(\intR  R(\beta_2 \beta_1 \beta , u) dx\right)^{\frac{1}{t}}||u||_{rt'}^{r}.
$$
Write $  R(\beta_2 \beta_1 \beta , u) = R\left(\beta_2 \beta_1 \beta ||u||^{\frac{N}{N-s}} , \frac{u}{||u||}\right) $. Then, choosing $ M > 0 $ sufficiently small such that $ \beta_2 \beta_1 \beta ||u||^{\frac{N}{N-s}} < \crit $, from Lemma \ref{l2}, we can find a constant $ K_1 = K_1(\beta , N, M, l, s) > 0 $ satisfying
$$
\intR |u|^r R(\beta , u)^l dx \leq K_1 ||u||_{rt'}^{r} .
$$
Setting $ t_0 = rt' > p $, we have the desired result. 
\end{proof}

\bo\label{obs3}
Under the hypothesis of Lemma \ref{l11}, from the embeddings \cite[Theorem 6.9]{[hitch]} $ t_0 > p $, there exists a constant $ K_2 = K_2(\beta , N, M, l, s) > 0 $ such that
$$
\intR |u|^r R(\beta , u)^l dx \leq K_2 ||u||^{r} .
$$
Moreover, since $ ||\cdot||\leq ||\cdot||_X $, we have
$$
\intR |u|^r R(\beta , u)^l dx \leq K_2 ||u||_{X}^{r} .
$$
\eo

\bo\label{obs6}
Observe that, the estimates obtained in Lemma \ref{l11} and Remark \ref{obs3} can be applied for an arbitrary, but fixed, $ u\in X\setminus\{0\} $ making $ \beta > 0 $ sufficiently small in order to apply Moser-Trudinger Lemma \ref{l2}.
\eo

\bl\label{l12}
There exists $ \rho >0 $ such that
\begin{equation}\label{eq14}
m_\beta = \inf\{ I(u) \ ; \ u\in X \ , \ ||u||=\beta \} > 0 \ , \ \forall \ \beta \in (0, \rho] 
\end{equation}
and
\begin{equation}\label{eq15}
n_\beta = \inf\{ I'(u)(u) \ ; \ u\in X \ , \ ||u||=\beta \} > 0 \ , \ \forall \ \beta \in (0, \rho]. 
\end{equation}
\el
\begin{proof}
Let $ u\in X\setminus \{0\} $, with $ ||u|| $ sufficiently small in order to apply Lemma \ref{l11}, and $ q>p $. Then, from (\ref{v2}), Remark \ref{obs3}, (\ref{eq5}) and Sobolev embeddings \cite[Theorem 6.9]{[hitch]}, we have
$$
I(u) \geq \dfrac{1}{p}||u||^p - \dfrac{K_0}{2p}||u||_{\frac{2N}{2N-1}p}^{2p}- \dfrac{\varepsilon}{p}||u||_{p}^{p}-K_2 ||u||^q \geq \dfrac{||u||^p}{p}[1-\varepsilon - C_1||u||^p - C_2 ||u||^{q-p}] .
$$
Hence, for $ \varepsilon > 0 $ and $ \rho > 0 $ sufficiently small, we obtain (\ref{eq14}). Similarly, from (\ref{eq7}), (\ref{v2}), Remark \ref{obs3} and Sobolev embeddings \cite[Theorem 6.9]{[hitch]}, follows
$$
I'(u)(u)=||u||^p +V_1(u)-V_2(u)-\intR f(u)u dx \geq ||u||^p [1-\varepsilon -C_3||u||^p -C_4||u||^{q-p}] .
$$
Therefore, taking $ \varepsilon, \rho > 0 $ sufficiently small, we get (\ref{eq15}).
\end{proof}

\bl\label{l13}
Let $ u\in X\setminus \{0\} $, $ t>0 $ and $ q>2p $. Then,
$$
\lim\limits_{t\RA 0} I(tu) = 0 \ \ , \ \ \sup\limits_{t>0}I(tu) < + \infty \mbox{ \ \ and \ \ } I(tu)\RA - \infty \mbox{ \ as \ } t\RA + \infty .
$$
\el
\begin{proof}
Let $ u\in X\setminus \{0\} $. First of all, from $ (f_4) $,
$$
I(tu)= \dfrac{t^p}{p}||u||^p + \dfrac{t^{2p}}{2p}V_0(u) - \intR F(tu) dx \leq \dfrac{t^p}{p}||u||^p + \dfrac{t^{2p}}{2p}V_0(u) - C_q t^q ||u||_{q}^{q}\RA - \infty ,
$$
as $ t\RA + \infty $. Now, from (\ref{eq5}) and Lemma \ref{l11}, for $ t>0 $ sufficiently small such that $ ||tu|| $ is under Lemma \ref{l11} conditions, we have
$$
\left| \intR F(tu) dx \right| \leq \dfrac{t^p}{p}||u||^p + K_1 t^q ||u||_{t_0}^{q} \RA 0 ,
$$
as $ t\RA 0 $. Hence, $ I(tu) \RA 0 $ as $ t\RA 0 $. Finally, once $ I\in C^{1}(X, \R) $, from the above two facts we conclude that $ \sup\limits_{t>0}I(tu) < + \infty $.
\end{proof}

Consider a sequence $ (\un) \subset X $ satisfying 
\begin{equation}\label{eq16}
\exists \ d > 0 \mbox{ \ s.t. \ } I(\un) < d \ , \ \forall \ n \in \N \mbox{ \ and \ } ||I'(\un)||_{X'}(1+||\un||_X) \RA 0 \mbox{ \ , as \ } n\RA + \infty .
\end{equation}

\bl\label{l14}
Let $ (\un) \subset X $ satisfying (\ref{eq16}). Then, $ (\un) $ is bounded in $ \Wsp $.
\el
\begin{proof}
From (\ref{eq16}) and $ (f_3) $, we have
$$
d + o(1) \geq I(\un) - \dfrac{1}{2p}I'(\un)(\un) \geq \dfrac{1}{2p}||\un||^p + \left(\dfrac{\theta}{2p}-1\right) \intR F(\un) dx \geq \dfrac{1}{2p}||\un||^p ,
$$
for all $ n\in \N $. Therefore, $ 2pd + o(1) \geq ||\un||^p $, for all $ n \in \N $, and the lemma follows.
\end{proof}

\bo\label{obs7}
\textbf{(1)} Observe that exchanging the condition $ I(\un)\leq d $, for all $ n\in \N $, for $ I(\un)\RA d >0 $ the above result remains valid.

\noindent \textbf{(2)} One can easily verify that the value $ c_{mp} $ satisfies $ 0<m_\rho \leq c_{mp} <+\infty $.

\noindent \textbf{(3)} Since $ I $ has the mountain pass geometry and $ c_{mp}>0 $, one can prove, as in \cite{[boer]}, that there exists a sequence $ (\un)\subset X $ such that 
\begin{equation}\label{eq17}
I(\un)\RA c_{mp} \mbox{ \ \ \  and \ \ \ } ||I'(\un)||_{X'}(1+||\un||_X)\RA 0 .
\end{equation}
Moreover, such sequence satisfies (\ref{eq16}).
\eo

\bl\label{l15}
Let $ (\un)\subset X $ satisfying (\ref{eq17}) and $ q> 2p $. Then, for some $ \rho_0 > 0 $ sufficiently small,
$$
\limsup\limits_{n} ||\un||^p < \rho_{0}^{p} .
$$
\el
\begin{proof}
From Lemma \ref{l14}, $ 2pc_{mp}+o(1) \geq ||\un||^p $, for all $ n \in \N $. Then, $ \limsup\limits_{n} ||\un||^p \leq 2pc_{mp} $. So, the natural step consists in finding a estimative for $ c_{mp} $. 

Consider the set $ \A = \{ u \in X\setminus \{0\} \ ; \ V_0(u) \leq 0 \} $. For each $ u\in X \setminus\{0\} $, $ t>0 $ and $ x \in \RN $, we define $ u_t(x)=t^2u(tx) $. Then, a directly computation give us
$$
V_0(u_t)=t^{4p-2N}V_0(u) - t^{4p-2N}\ln t ||u||_{p}^{2p} \RA - \infty ,
$$
as $ t\RA + \infty $, since $ 4p-2N=2p(2-s)>0 $. Hence, $ \A \neq \emptyset $. 

Moreover, from the embeddings \cite[Theorem 6.9]{[hitch]}, there exists $ C>0 $ such that $ ||u||\geq C||u||_q $. So, it makes sense to define
$$
S_q(v) = \dfrac{||v||}{||v||_q} \mbox{ \ \ \ and \ \ \ } S_q = \inf\limits_{v \in \A} S_q(v) \geq \inf\limits_{v \neq 0} S_q(v) > 0 .
$$

Now, from Lemma \ref{l13}, for $ v\in \A $ and $ T> 0 $ sufficiently large, $ I(Tv)<0 $. Defining $ \gamma : [0, 1] \RA X $ by $ \gamma(t)=tTv $, we have that $ \gamma \in \Gamma $ and
$$
c_{mp}\leq \max\limits_{0\leq t \leq 1} I(\gamma(t)) = \max\limits_{0\leq t \leq 1} I(tTv) \leq \max\limits_{t\geq 0} I(tv) .
$$
Consequently, for $ \psi \in \A $,
$$
c_{mp}\leq \max\limits_{t\geq 0} I(t\psi) \leq \max\limits_{t\geq 0} \left\{ \dfrac{t^p}{2p}||\psi||^p - C_q t^q ||\psi||_{q}^{q}\right\} \leq \left(\dfrac{1}{p}-\dfrac{1}{q}\right)\dfrac{S_q(\psi)^{\frac{pq}{q-p}}}{(qC_q)^{\frac{p}{q-p}}}. 
$$
Taking the infimum over $ \psi \in \A $, we obtain
$$
\limsup\limits_{n}||\un||^p \leq \dfrac{2(q-p)}{q}\dfrac{S_{q}^{\frac{pq}{q-p}}}{(qC_q)^{\frac{p}{q-p}}} \leq \rho_{0}^{p} ,
$$
for $ C_q >0$ sufficiently large. 
\end{proof}

\bo\label{obs8}
From Lemma \ref{l15}, one can see that, taking $ C_q > 0 $ sufficiently large, we can make $ \rho_0 $ sufficiently small such that every result concerning the exponential term is valid for sequences satisfying $ I(\un) \RA c_{mp} $ or $ I(\un) \leq c_{mp} $.
\eo

\section{Proof of Theorem \ref{t1}}

In the present section, we finish the proof of Theorem \ref{t1}. We start proving a key proposition that provides us with nontrivial critical points for $ I $ in $ X $. Then, it remains just to gather all the results inside the theorem. 

\bp\label{p3}
Let $ (f_1)-(f_4) $, $ q>2p $ and $ (\un)\subset X $ satisfying (\ref{eq17}). Then, up to a subsequence, only one between the two alternatives is valid
\begin{itemize}
\item[(a)] $ ||\un||\RA 0 $ and $ I(\un)\RA 0 $, as $ n \RA + \infty $.

\item[(b)] There exists points $ \yn \in \mathbb{Z}^N $ such that $ \until = \yn \ast \un \RA u $ in $ X $, for a nontrivial critical point $ u\in X $ of $ I $.
\end{itemize}
\ep
\begin{proof}
From Lemma \ref{l14}, $ (\un)\subset \Wsp $ is bounded. Suppose that (a) does not happen. 

\noindent \textbf{Claim 1:} $ \liminf\limits_{n\RA + \infty} \sup\limits_{y\in \mathbb{Z}^N} \ds_{B_2(y)} |\un(x)|^p dx > 0 $.

Lets suppose the contrary. Then, from and easy adaptation Lion's Lemma \cite[Lemma 2.4]{[yu]}, $ \un \RA 0 $ in $ \Lomega $, for all $ \omega > p $. Thus, since $ \frac{2N}{2N-1}p > p $, from (\ref{v2}), $ V_2(\un) \RA 0 $. Moreover, from (\ref{eq7}) and Remark \ref{obs3}, we have
$$
\left| \intR f(\un)\un dx \right| \leq \varepsilon ||\un||_{p}^{p}+C_1 ||\un||_{qt_0}^{q} \leq \varepsilon C_2 + C_1 ||\un||_{qt_0}^{q} \RA 0,
$$ 
as $ \varepsilon \RA 0 $ and $ n \RA + \infty $. Consequently, 
$$
||\un||^p +V_1(\un) = I'(\un)(\un) + V_2(\un) +  \intR f(\un)\un dx \RA 0 ,
$$
as $ n \RA + \infty $. So, from the non-negativeness, we have $ ||\un||\RA 0 $ and $ V_1(\un) \RA 0 $ and, from the embeddings \cite[Theorem 6.9]{[hitch]}, $ ||\un||_p \RA 0 $ and  $ ||\un||_{qt_0}\RA 0 $. Finally, from (\ref{eq5}), Remark \ref{obs3} and   the embeddings \cite[Theorem 6.9]{[hitch]}, we conclude that $ \intR F(\un) dx \RA 0 $. Hence, $ I(\un) \RA 0 $, which is a contradiction, proving the claim.

As a consequence, one can easily obtain, up to subsequence, $ (\yn) \subset \mathbb{Z}^N $ such that $ (\yn \ast \un) \subset X $ and $ \until = \yn \ast \un \CF u $ in $ \Wsp\setminus\{0\} $. So, without loss of generality, we can assume that $ \until(x) \RA u(x) $ a.e. in $ \RN $ and, from \cite[Theorem 6.9]{[hitch]}, $ (\until) $ is bounded in $ \Lomega $ for all $ \omega \geq p $.

Now, observe that, since $ \frac{2N}{2N-1}p , qt_0 > p $, we have
\begin{align*}
V_1(\until)=V_1(\un)& = I'(\un)(\un) + V_2(\un) + \intR f(\un) \un dx -||\un||^p \\
& \leq K_0||\un||_{\frac{2N}{2N-1}p}^{2p}+ \varepsilon C_2 C_1 ||\un||_{qt_0}^{q}+ o(1) \leq C_3 + o(1) . 
\end{align*}
That is, $ \sup\limits_{n}V_1(\until) < + \infty $. So, from Proposition \ref{p2}, $ (||\until||_\ast) $ is bounded. Since $ (\until) $ is already bounded in $ \Wsp $, $ (\until) $ is bounded in $ X $. Hence, from the reflexiveness of $ X $, passing to a subsequence if necessary, $ u\in X $ and $ \until \CF u $ in $ X $. Moreover, from Proposition \ref{p1}, $ \until \RA u $ in $ \Lomega $, for all $ \omega \geq p $.

\noindent \textbf{Claim 2:} $ I'(\until)(\until - u) \RA 0 $, as $ n\RA + \infty $.

First of all, observe that, by a change of variables, $ I'(\until)(\until - u) = I'(\un)(\un - (-\yn)\ast u) $. Thus,
\begin{equation}\label{eq18}
|I'(\until)(\until - u)| = |I'(\un)(\un - (-\yn)\ast u)|\leq ||I'(\un)||_{X'}(||\un||_X + ||(-\yn)\ast u||_X).
\end{equation}
Then, similarly as in \cite{[boer]}, we first seek for an useful inequality for $ ||(-\yn) \ast u||_X $. If $|\yn|\RA +\infty$, then, for $ x\in \RN $,
$$
\ln(1+|x-\yn|)-\ln(1+|\yn|)=\ln \left(\dfrac{1+|x-\yn|}{1+|\yn|}\right) \RA 0 , n\RA + \infty .
$$
Therefore, there exists $ C_4 > 0 $ such that $ \ln(1+|x-\yn|)\geq C_4 \ln(1+|\yn|) $. \\
Now, suppose that $ (\yn)\subset \mathbb{Z}^N $ converges to $ y_0 \in \mathbb{Z}^N $. Then, up to a subsequence, $ \yn \equiv y_0 $. Let $ y_0 \neq 0 $ and consider $ r $ the line passing through the origin and $ y_0 $. Then, define $ \Omega_0 $ as the open connected region between $ r $ and one of the axis, such that the angle between $ r $ and the axis is $ \leq \frac{\pi}{2} $. Taking $ \delta > 0 $ such that $ \delta < |y_0| $, set $ \Omega = \Omega_0 \cap B_\delta $. So, for $ x\in \Omega $, we have that $ |x-y_0| > |y_0| $. Therefore, by the Mean Value Theorem, there exists $ x_\delta \in \Omega $, satisfying
\begin{align*}
||\un||_{\ast}^{2} & \geq \int\limits_{\Omega} \ln(1+|x-\yn|)|\until (x)|^p dx \\
& = |\Omega| \ln(1+|x_\delta -\yn|)|\until (x_\delta )|^p \\
& = C_4 \ln(1+|x_\delta -y_0|) \geq C_4 \ln(1+|y_0|) = C_4 \ln(1+|\yn|) ,
\end{align*}
for $ C_4 > 0 $. If $ y_0 =0 $, $ \until = \un $ and the result follows immediately from (\ref{eq18}). So, in any of the cases, there exists $ C_4 >0 $ such that
$$
||\un||_{\ast}^{p}= \intR \ln(1+|x -\yn|)|\until (x)|^p dx \geq C_4 \ln(1+ |\yn|) \ , \forall \ n \in \N .
$$
Now, we have
$$
||\until||_{\ast}^{p} = \intR \ln (1+|x+\yn|) |\un(x)|^p dx \leq ||\un||_{\ast}^{p} + \ln (1+|\yn|) ||\un||_{p}^{p} .
$$
From this, since every norm is weakly lower semicontinuos, $ ||\cdot||_p $ is $ \mathbb{Z}^N $-invariant and $ \until \rightharpoonup u $ in $ X $, follows that
\begin{align*}
||(-\yn) \ast u||_{\ast}^{p}&=\intR \ln(1+|x -\yn|)|u (x)|^p dx \\
& \leq ||\until||_{X}^{p} + \ln(1+|\yn|)||\un||_{p}^{p} \\
& \leq ||\un||^{p}  + ||\un||_{\ast}^{p} (1 + C_5||\un||_{p}^{p}) \\
& \leq ||\un||^{p}  + C_6 ||\un||_{\ast}^{p} \leq C_7 ||\un||_{X}^{p}
\end{align*}
for $ n\in \N $ and $ C_7 >0 $. Consequently, there exists a constant $ C_8 >0 $ such that, after passing to a subsequence, we have, for all $ n\in \N $,
\beq\label{eq29}
||(-\yn) \ast u||_X^{p} =||u||^p + ||(-\yn) \ast u||_{\ast}^{p} \leq ||\un||^p +C_7 ||\un||_{X}^{p} \leq C_8 ||\un||_{X}^{p} .
\eeq
Therefore, from (\ref{eq18}) and (\ref{eq29}),
$$
|I'(\until)(\until - u)| \leq (1+C_{8}^{\frac{1}{p}})||I'(\un)||_{X'}||\un||_X \RA 0 ,
$$
finishing the claim. 

\noindent \textbf{Claim 3:} $ \intR f(\until)(\until - u) dx \RA 0 $, as $ n \RA + \infty $.

Since $ ||\cdot|| $ is $ \mathbb{Z}^N $-invariant, Lemma \ref{l2}, Lemma \ref{l11} and Remark \ref{obs3} remains valid for $ (\until) $. Moreover, $ \frac{1}{2p}+ \frac{2p-1}{2p}=1 $, $ \frac{2p}{2p-1}(q-1)> p $ and $ \omega_0 = \frac{2p}{2p-1}(q-1)t_0 > p $. Thus, from the boundedness of $ (\until) $ in $ L^{\omega_0}(\RN) $ and in $ \Lp $, Remark \ref{obs3}, Proposition \ref{p1}, (\ref{eq7}) and Hölder inequanlity, we have
\begin{align*}
\left| \intR f(\until)(\until - u) dx \right| & \leq \intR |\until|^{p-1}|\until - u| dx + b_1 \intR |\until|^{q-1}|\until - u| R(\alpha , \until) dx \\
& \leq C||\until - u||_p + b_1||\until -u ||_{2p} (\textstyle{\int_{\RN}} R(\alpha , \until)^{\frac{2p}{2p-1}}|\until|^{\frac{2p}{2p-1}(q-1)} dx )^{\frac{2p-1}{2p}} \\
& \leq  C||\until - u||_p + b_1 K_1 ||\until -u ||_{2p} ||\until||_{\omega_0}^{q-1}\RA 0 ,
\end{align*} 
as $ n\RA + \infty $, proving the claim. 

Moreover, we observe that

\noindent \textbf{(i)} From HLS, Hölder inequality and Proposition \ref{p1}, 
$$
|V_2'(\until)(\until -u)|\leq K_0 ||\until||_{\frac{2N}{2N-1}p}^{2p-1}||\until - u||_{\frac{2N}{2N-1}p}\RA 0 .
$$

\noindent \textbf{(ii)} $ \left| \intR |\until|^{p-2}\until (\until - u) dx \right| \leq ||\until||_{p}^{p-1}||\until - u||_p \RA 0 $. 

\noindent \textbf{(iii)} We recall a standard result: there exists a constant $ D_1>0 $, depending only on $ p $, such that
\begin{equation}\label{eq21}
|a-b|^p \leq D_1 (|a|^{p-2}a-|b|^{p-2}b)(a-b) \ , \ \forall \ a, b\in \R \ , \ \forall \ p \geq 2 .
\end{equation}
So, 
\begin{align*}
V_1'(\until)(\until - u) & = \intR \intR \ln(1+|x-y|)|\until(x)|^p |\until(y)|^{p-2}\until(y)(\until(y)-u(y))dxdy \\
& \geq C_9 \intR \intR \ln(1+|x-y|)|\until(x)|^p |\until(y)-u(y)|^{p} dxdy \\
& + \intR \intR \ln(1+|x-y|)|\until(x)|^p |u(y)|^{p-2}u(y)(\until(y)-u(y))dxdy \\
& = C_9 A_1 + B_1 .
\end{align*}
Note that $ A_1 \geq 0 $ and, from Lemma \ref{l9}, $ B_1\RA 0 $, as $ n\RA + \infty $. Hence, $ V_1'(\until)(\until - u)\RA 0 $.

Hence, from Claims 2 and 3 and items (i)-(iii), we conclude that
\begin{align*}
o(1) = I'(\until)(\until - u) & = \tilde{A}(\until)(\until - u) + V_0 '(\until(\until -u) - \intR f(\until)(\until -u) dx \\
& \geq \tilde{A}(\until)(\until -u) + o(1) .
\end{align*}
That is, $ \tilde{A}(\until)(\until -u)\RA 0 $. So, once $ \tilde{A} $ has the $ (S) $ property, $ \until \RA u $ in $ \Wsp $. Therefore, we can also conclude that $ A_1 \RA 0 $ and, from Proposition \ref{p2}, $ ||\until - u ||_\ast \RA 0 $, proving that $ \until \RA u $ in $ X $. 

Finally, remains to show that $ u $ is a critical point of $ I $. Let $ v\in X $. So, as we did above, is possible to find $ C_{10} > 0 $ such that $ ||(-\yn)\ast v||_X \leq C_{10} ||\un||_X $. Thus,
$$
|I'(u)(v)|=\lim |I'(\until)(v)| = \lim |I'(\un)((-\yn)\ast v)| \leq C_{10} \lim ||I'(\un)||_{X'}||\un||_X = 0 . 
$$
Therefore, $ u $ is a nontrivial critical point for $ I $ in $ X $.
\end{proof}

\begin{proof}[Proof of Theorem \ref{t1}]
\textbf{(i)} From Lemma \ref{l12} and Proposition \ref{p3} there exists a nontrivial critical point of $ I $, $ u_0 \in X $, such that $ I(u_0)=c_{mp} $. 

\textbf{(ii)} We start defining the set $ \K = \{v \in X\setminus \{0\} \ ; \ I'(v)=0\} $. Since $ u_0\in \K $, $ \K \neq \emptyset $. Thus, we can consider a sequence $ (\un)\subset \K $ satisfying $ I(\un)\RA c_g = \inf\limits_{v\in \K } I'(v) $.

Observe that $ c_g \in [-\infty , c_{mp}] $. If $ c_q = c_{mp} $ the proof is finished. Otherwise, if $ c_g < c_{mp} $, considering a subsequence if necessary, we can suppose that $ I(\un)\leq c_{mp} $, for all $ n\in \N $ and, from the definition of $ \K $, we see that $ (\un) $ satisfies $ ||I'(\un)||_{X'}(1+||\un||_X)\RA 0 $. Moreover, since $ I'(\un)(\un) = 0 $, for all $ n\in \N $, from (\ref{eq15}), $ ||\un|| > \rho $, for all $ n \in \N $. Therefore, from Proposition \ref{p3} there exists $ (\yn)\subset \mathbb{Z}^N $ such that $ \until \RA u $ in $ X $, for a nontrivial critical point $ u $ of $ I $ in $ X $. Consequently, $
I'(u)=\lim I'(\until)=\lim I'(\un)=0 
$
and we conclude that $ u\in \K $ and 
$$
I(u)=\lim I(\until) = \lim I(\un) = c_g .
$$
Particularly, $ c_g > -\infty $.
\end{proof}

\section{Proof of Theorem \ref{t2}}

In this section we will provide the proof of the multiplicity result stated in Theorem \ref{t2}. In order to do so, we will need to verify some results concerned with the genus theory, denoted by $ \gamma $ and whose definition, given over $ \A = \{ A\subset X \ ; \ A \mbox{ is symmetric and closed} \} $ (with respect to continuity in $ X $), and basic properties can be found in \cite{[genus]} Chapter II.5. We start this section, verifying some properties of an important auxiliary function, namely, $ \varphi_u : \R \RA \R $, given by $ \varphi_u(t) = I(tu) $, for all $ u\in X\setminus\{0\} $ and $ t \in \R $.

Once the results are done similarly as in \cite{[boer]}, in order to make the paper concise, we will only sketch the proofs here (see also \cite{[6]}).  

\bl\label{l110}
\textbf{(a)} Let $ u\in X\setminus \{0\} $. Then, $ \varphi_u $ is even and there exists a unique $ t_u \in (0, +\infty) $ such that $ \varphi_u '(t)>0 $, for all $ t\in (0, t_u) $, and $ \varphi_u '(t)<0 $, for all $ t\in (t_u , \infty) $. Moreover, $ \varphi_u (t) \RA -\infty $, as $ t \RA +\infty $.

\noindent \textbf{(b)} Let $ u\in X\setminus \{0\} $. Then, there exists a unique $ t_{u}'\in (0, + \infty )$ such that $ \varphi_u(t)>0 $, for $ t\in (0, t_{u}') $, and $ \varphi_u(t)<0 $, for $ t\in (t_{u}', + \infty) $. Moreover, $ t_u $ given by item (a) is a global maximum for $ \varphi_u $.

\noindent \textbf{(c)} For each $ u\in X\setminus \{0\} $, the map $ u \mapsto t_u ' $ is continuous.
\el
\begin{proof}
We prove item (a). Since $ f $ is odd, $ I $ is even and, consequently, $ \varphi_u $ is even as well. 

\noindent \textbf{(i)} For $ t > 0 $ sufficiently small and $ \alpha > 0 $, for (\ref{eq7}), we have
$$
\varphi_u '(t) \geq t^{p-1}||u||^p[1- C_2t^{2(p-1)} ||u||^{p} -\varepsilon -C_4t^{q-p}||u||^{q-p}] .
$$
Thus, $ \varphi_u '(t) > 0 $ for $ t, \varepsilon > 0 $ sufficiently small. 

\noindent \textbf{(ii)} From $ (f_3) $ and $ (f_4 ') $, once $ q>2p $, 
$$
\varphi_u '(t)  \leq t^{p-1}||u||^p + t^{2p-1} V_1(u) - C_3 t^{q-1}||u||_{q}^{q} \RA - \infty \mbox{ \ , as \ } t\RA + \infty .
$$

Hence, from (i)-(ii), since $ I $ is $ C^1 $, there exists $ t_u \in (0, + \infty) $ such that $ \varphi_u '(t_u) =0 $, which is unique by $ (f_5) $. 

Item (b) follows as a consequence of item (a) and item (c) as a consequence of predecessors. 
\end{proof}

Now, we define the following sets 
$$
K_c = \{ u\in X \ ; \ I'(u)=0 , I(u)=c \} \ , \  c\in (0, + \infty)
$$
and
$$
A_{c, \rho} \{ u\in X \ ; \ ||u-v|| \leq \rho \mbox{ \ , for some } v\in K_c\} \ , \  c\in (0, + \infty).
$$
It is easy to verify that the sets $ K_c $ and $ A_{c, \rho} $ are symmetric, closed (with respect to $ X $) and invariant under $ \mathbb{Z}^N $ translations, i.e, if $ u\in K_c , A_{c, \rho} $, then $ z\ast \in K_c , A_{c, \rho} $, for all $ z\in \mathbb{Z}^N $.

Next, we fix a continuous map $ \beta: \Lp \setminus \{0\} \RA \RN $ that is equivariant under $ \mathbb{Z}^N $ translations, that is, $ \beta(x \ast u) = x+\beta(u) $, for $ x\in \mathbb{Z}^N $ and $ u\in \Lp \setminus \{0\} $.  We also require that $ \beta(-u)=\beta(u) $. Such map is called a \textit{generalized barycenter map} and an example can be constructed as in \cite{[bar]}. Hence, we can define 
$$
\tilde{K}_c = \{ u\in K_c \ ; \ \beta(u)\in [-4, 4]^N \} ,
$$
which are clearly symmetric sets. Moreover, before given our first result, we need to recall the Gauss bracket $ [ \cdot ] : \R \RA \mathbb{Z} $, given by $ [s]=\max\{n\in \mathbb{Z} \ ; \ n \leq s \} $ (see \cite[Chapter 3]{[n2]}), which naturally induce a map from $ \RN $ onto $ \mathbb{Z}^N $, as follows
$$
[(x_1, x_2, ..., x_N)] = ([x_1], [x_2], [x_N])\in \mathbb{Z}^N \ , \ \forall \ (x_1, x_2, ..., x_N)\in \RN .
$$
We recall some properties of the Gauss bracket, that are needed inside of the proof of our results.

\bl\label{l18}
Let $ [ \cdot ] : \R \RA \mathbb{Z} $, given by $ [s]=\max\{n\in \mathbb{Z} \ ; \ n \leq s \} $. Then,

\noindent \textbf{(i)} $ 0 \leq s - [s] < 1 $, for all $ s\in \R $.

\noindent \textbf{(ii)} if $ z\in \mathbb{Z} $ and $ s\in \R $, then $ [z+s]=z+[s] $.

\noindent \textbf{(iii)} Let $ s\in \R $ such that $ s - [s] \geq \dfrac{1}{2} $, then $ s - \dfrac{1}{2} - \left[ s - \dfrac{1}{2}\right] < \dfrac{1}{2} $.

\noindent \textbf{(iv)} $ 0 \leq s - [s] < 1 $, for all $ s\in \R $.

\noindent \textbf{(v)} Let $ s\in \R $. Then, $ s- \left[s - \dfrac{1}{2}\right]< 1 $.
\el

\bp\label{p19}
Let $ c>0 $. Then, there exists $ \rho_0 = \rho_0(c) > 0 $ such that $ \gamma(A_{c, \rho}) < \infty $, for all $ \rho \in (0, \rho_0) $.
\ep
\begin{proof}
The proof of this proposition can be done similarly as \cite[Proposition 4.1]{[boer]}, with minor changes. We only highlight how to construct the sets $ L_i $ inside the refereed proof, once here we need $ 2^N $ sets to ``cover'' $ \Lp \setminus \{0\} $.
 
Let $ \{e_i\}_{i=1}^{N} $ be the canonical base of $ \RN $. Consider $ a_{i} = \frac{1}{2}e_i $, $ a_{ij} = \frac{1}{2}e_i + \frac{1}{2}e_j $, $ a_{ijk}= \frac{1}{2}(e_i+e_j+e_k) $, for $ 1\leq i , j, k \leq N $, and successively until $ a_{12\cdots N}=\frac{1}{2}\sum\limits_{i=1}^{N}e_i $. Observe that $ a_{12}=a_{21} $, so, excluding the repeating cases, we define $ 2^N-1 $ sets as $ L_{i} = a_{i} \ast L_1 \subset \Lp \setminus \{0\} $, $ L_{ij} = a_{ij} \ast L_1 \subset \Lp \setminus \{0\} $, ..., $ L_{12\cdots N} = a_{12\cdots N} \ast L_1 \subset \Lp \setminus \{0\} $. Once this construction is well understood, in order to simplify the notation, we simply denote $ a_i $ and $ L_i $, for $ 2 \leq i \leq 2^N $. 

It is clear that $ \Lp \setminus \{0\} \subset \bigcup\limits_{i=1}^{2^N}L_i $. 
\end{proof}

For the next results we will need the definition and some basic properties of relative genus. So, for convenience of the reader, we will include it here.

\bd\label{d11}
Let $ D, Y \in \A $ with $ D \subset Y $. We say that $ U, V\in \A $ is a covering of $ Y $ relative to $ D $ if is satisfies

\noindent \textbf{(i)} $ Y\subset U \cup V $ and $ D\subset U $;

\noindent \textbf{(ii)} there exists an even continuous (in $ X $) function $ \chi : U \RA D $, such that $ \chi(u)=u $, for all $ u\in D $. 

If $ U, V\in \A $ is a covering of $ Y $ relative to $ D $, then the genus of this covering is $ \gamma(V)=k $. 
\ed 

\bd\label{d12}
Let $ D, Y \in \A $ with $ D \subset Y $. We define the Krasnoselskii's Genus of $ Y $ relative to $ D $, denoted by $ \gamma_D(Y) $, as

\noindent \textbf{(i)} There exists a covering for $ Y $ relative to $ D $ and, in this case, $ \gamma_D(Y)=k $, where $ k $ is the lowest genus of this coverings. 

\noindent \textbf{(ii)} If we cannot find any such covering of $ Y $ relative to $ D $, we set $ \gamma_D(Y)=+\infty $.
\ed

In the following, we list some useful properties of relative genus that are needed to guarantee that the results are valid. 

\bl\label{l19}
\textbf{(i)} Let $ D \subset \A $. Then, $ \gamma_D(D)=0 $.

\noindent \textbf{(ii)} Let $ D, Y, Z \in \A $ satisfying $ D\subset Y $ and $ D\subset Z $. If there exists a function $ \varphi : Y \RA Z $, even and continuous (in $ X $), such that $ \varphi(u)=u $, for all $ u\in D $, then $ \gamma_D(Y)\leq \gamma_D(Z) $.

\noindent \textbf{(iii)} Let $ D\subset Y\subset Z \in \A $. Then, $ \gamma_D(Y)\leq \gamma_D(Z) $.

\noindent \textbf{(iv)} Let $ D, Y, Z \in \A $ satisfying $ D\subset Y $. Then, $ \gamma_D(Y\cup Z)\leq \gamma_D(Y)+\gamma(Z) $.
\el
\begin{proof}
For item (i), take $ U=D , V= \emptyset $ and $ \chi = id $, in the definition of relative genus. Proofs for itens (ii) and (iv) can be found, for example, in \cite{[6]}. Finally, item (iii) is an immediate consequence of item (ii).
\end{proof}

For the next results, we define the sets
$$
I^{c}=\{u\in X \ ; \ I(u)\leq c \} \mbox{ \ , for \ } c\in \R \mbox{ \ and \ } D=I^0 ,
$$
and the values
$$
c_k = \inf\{c \geq 0 \ ; \ \gamma_D(I^c)\geq k \} \mbox{ \ , \ } \forall \ n\in \N .
$$

\bo\label{obs14}
\textbf{(1)} Since $ I $ is unbounded from bellow, $ D\neq \emptyset $. 

\noindent \textbf{(2)} Let $ c_1, c_2 \in \R $ with $ c_1 > c_2 $. Then, if $ u\in I^{c_2} $, $ I(u)\leq c_2 < c_1 $, so $ u\in I^{c_2} $. That is, if $ c_1 > c_2 $, then $ I^{c_2}\subset I^{c_1} $.

\noindent \textbf{(3)} If $ c_1 > c_2 \geq 0 $, then $ D\subset I^{c_2}\subset I^{c_1} $. Consequently, $ \gamma_D(I^{c_2})\leq \gamma_D(I^{c_1}) $.

\noindent \textbf{(4)} For $ \varepsilon >0 $, $ \gamma_D(I^{c_k + \varepsilon})\geq k $ and $ \gamma_D(I^{c_k - \varepsilon})<k $, for every $ k \in \N $.

\noindent \textbf{(5)} $\inf\limits_{u\in X \setminus \{0\}} \sup\limits_{t \in \R} I(tu) = \inf\limits_{u\in X \setminus \{0\}} \sup\limits_{t >0 } I(tu) < + \infty .$
\eo

Consider the Nehari's manifold for $ I $, defined by
\begin{equation}\label{VN}
\Nn = \{ u\in X \setminus \{0\} \ ; \ I'(u)(u) = 0 \}.
\end{equation}

\bl\label{l113}
Let $ \Nn $ as in (\ref{VN}). Then, 
$$
\inf\limits_{\Nn} I = \inf\limits_{u \in X \setminus \{0\}}\sup\limits_{t > 0} I(tu) .
$$
\el
\begin{proof}
From Lemma \ref{l110} and the chain rule, for $ u\in X \setminus \{0\} $, we have $ t_u u \in \Nn $ and $ \sup\limits_{t > 0} I(tu) = I(t_u u) $. Moreover, for $ u\in \Nn $, once again from Lemma \ref{l110}, $ \sup\limits_{t > 0} I(tu) = I(u) $. Therefore, the result follows.
\end{proof}

\bl\label{l114}
We have $\inf\limits_{X \setminus \{0\}} \sup\limits_{t\in \R} I(tu) > 0 $. 
\el
\begin{proof}
From Lemmas \ref{l113} and \ref{l12} and direct calculations, one can obtain the result.  
\end{proof}

In order to prove the next results, we need to introduce the following sets
$$
\Nn^{+} = \{ u\in X \ ; \ I'(u)(u) > 0 \} \mbox{ \ \ and \ \ } \Nn^{-} = \{ u\in X \ ; \ I'(u)(u) < 0 \}.
$$
Note that $ X = \{0\} \cupdot \Nn^{+} \cupdot \Nn \cupdot \Nn^{-} $. Moreover, using the definition of a set's boundary and Lemmas \ref{l12} and \ref{l110}, it is possible to verify that $ \partial \Nn^{-} = \Nn $ and $ \partial \Nn^{+} = \{0\} \cup \Nn $.

\bp\label{p110}
We have $ c_1 = \inf\limits_{N} I = \inf\limits_{u\in X \setminus \{0\}}\sup\limits_{t > 0}I(tu) >0 $.
\ep
\begin{proof}
Observe that, from Lemma \ref{l113}, remains to prove that $ c_1 = \inf\limits_{N} I > 0 $. From Lemma \ref{l19}-(i), $ c_1 >0 $.

\noindent \textbf{Claim 1:} $ c_1 \geq \inf\limits_{\Nn} I $.

Suppose, by contradiction, that $ c_1 < \inf\limits_{\Nn} I  $. Choose $ c\in ( c_1, \inf\limits_{\Nn} I) $. Define the function $ F: I^c \RA X $ by
$$
F(u) = \left\{ \begin{array}{ll}
0 \mbox{ \ , if \ } u\in \{0\} \cup \Nn^{+} \\
\max\{ 1, t_u '\}u \mbox{ \ , if \ } u\in \Nn^{-}
\end{array}
\right. .
$$
We see that $ F $ is well-defined, continuous, odd and $ F\big\vert_{D} = id $ Consider $ U= I^c $ and $ V = \emptyset $. Since $ I $ is $ C^1 $ and odd, $ U $ is closed and  symmetric. Hence, $ \gamma_D (I^c ) =0 $. But it gives a contradiction, since $ 1 \leq \gamma_D(I^{c_1}) \leq \gamma_D(I^c) = 0 $. Consequently, $ c_1 \geq \inf\limits_{\Nn} I $. 

\noindent \textbf{Claim 2:} $ c_1 \leq \inf\limits_{\Nn} I $.

For $ u_0 \in X\setminus\{0\} $, without loss of generality, we can assume $ ||u_0|| =1 $. Set $ d = \sup\limits_{t>0} I(tu_0) $. Lets prove that $ \gamma_D(I^d ) \geq 1 $. Note that, if $ u\in B= \{tu_0 \ ; \ t > 0\} $, then there exists $ t_0 > 0 $ such that $ u= t_0 u_0 $. One can see that $D \subset B \cup D \subset I^d $. 

Thus, from Lemma \ref{l19}-(iii), $ \gamma_D(B \cup D) \leq \gamma_D(I^d) $. In this point of view, we work to prove that $ \gamma_D(B \cup D) \geq 1 $. 

Suppose that  $ \gamma_D(B \cup D) =0 $. Then, by definition, $ U= B\cup D $ and $ V = \emptyset $, once only the empty set has null genus. Moreover, there exists a function, continuous and odd, $ \chi : B \cup D \RA D $ such that $ \chi(u)=u $, for all $ u\in D $. 

Define $ g: (0, + \infty) \RA (0, + \infty) $ by $ g(t) = ||\chi(tu_0)|| $. Note that $ g $ is continuous. Hence, from Lemmas \ref{l12} and \ref{l110} and the Intermediate Value Theorem, we reach a contradiction. Therefore, claim 2 is valid and, combined with claim 1, we have the proposition.  
\end{proof}

\bo\label{obs15}
\textbf{(1)} We can provide an equivalent definition to the the function $ F: I^c \RA X $ by
$$
F(u) = \left\{ \begin{array}{ll}
0 \mbox{ \ , if \ } u\in \{0\} \cup \Nn^{+} \\
\sigma(u) u \mbox{ \ , if \ } u\in \Nn^{-}
\end{array}
\right. ,
$$
where $ \sigma : \Nn^{-} \RA  [1, + \infty) $ is given by $ \sigma(u) = \inf\{t \geq 1 \ ; \ \varphi_u(t)=I(tu) \leq 0 \} $. 

\noindent \textbf{(2)} Let $ i, j\in \N $ with $ i>j $. Then, 
$$
\{c\in (0, + \infty) \ ; \ \gamma_D(I^c) \geq i \} \subset \{c\in (0, + \infty) \ ; \ \gamma_D(I^c) \geq j \} . 
$$
Therefore, $ c_j \leq c_i $. 
\eo

In the following, we will consider $ W $ as a $ k $-dimensional subspace of $ X $. As a consequence of $ (f_1 ') $, $ (f_4 ') $, $ V_2(u)\geq 0 $ and $ V_1(u)\leq 2||u||_{X}^{2p} $, once can see that 
\begin{equation}\label{l116}
I(u)\RA - \infty \mbox{, as \ }  ||u||_X \RA + \infty \mbox{, and \ } \sup\limits_{u\in W} I(u) < + \infty .
\end{equation}

\bc\label{l117}
There exists $ R>0 $ such that $ \{u\in W \ ; \ ||u|| \geq R \} \subset D $. 
\ec

For the next results, consider $ \rho > 0  $ given by Lemma \ref{l12} and $ \chi : X \RA X $ a continuous and even function, such that $ \chi(u)=u $, for all $ u\in D $. Define the sets
$$
\Oo_\chi = \{ u\in W \ ; \ ||\chi(u)|| < \rho \} .
$$ 
For convenience, we set $ \chi $ as a function having these properties, except we say otherwise. 

\bo\label{obs16}
One can easily verify that the sets $ \Oo_\chi $ have the following properties:

\noindent \textbf{(1)} $ \Oo_\chi $ is a neighbourhood of zero in $ W $. 

\noindent \textbf{(2)} $ \Oo_\chi $ is bounded and symmetric.

\noindent \textbf{(3)}  If $ u\in \partial_W \Oo_\chi $, then $ ||\chi(u)||=\rho $.
\eo

\bp\label{p111}
Let $ k\in \N $. Then $ c_k < + \infty $.
\ep
\begin{proof}
Set $ \alpha = \sup\limits_{u\in W} I(u) $. From (\ref{l116}), $ \alpha < + \infty $. Also, by definition of $ I^\alpha $, we see that $ W \subset I^\alpha $. 

Suppose that $ \gamma_D (I^\alpha) < k $. Then, from Definition \ref{d12}, a corollary of Tietze's Theorem, Remark \ref{obs16} and a genus property, $ \gamma(\partial_W \Oo_\chi ) = k $. 

Let $ u\in \partial_W \Oo_\chi  $. From Remark \ref{obs16}, $ ||\chi(u)||= \rho $. By Lemma \ref{l12}, $ I(\chi(u))\geq m_\rho > 0 $. Thus, $ \chi(u)\not\in D $ and $ \partial_W \Oo_\chi \cap U = \emptyset $. 

On the other hand, $$ \partial_W \Oo_\chi \subset W \subset I^{\alpha} \subset U \cup V \ \ \Longrightarrow \ \ \partial_W \Oo_\chi \subset V . $$
Then, we have $ k = \gamma(\partial_W \Oo_\chi) \leq \gamma(V) \leq k-1 $, which is a contradiction. 

Therefore, from the definition of values $ c_k $ and Lemma \ref{l116}, $ c_k \leq \alpha < + \infty $. 
\end{proof}

Before we prove that the values $ c_k $ are critical values of $ I $, we will provide a deformation lemma. We start defining the sets $ S = X \setminus A_{c, \rho} $, 
$$
S_\delta = \{ u \in X \ ; \ ||u-v||_X \leq \delta \mbox{ \ , for some \ } v\in S \} 
$$
and 
$$
\tilde{S}_\delta = \{ u \in X \ ; \ ||u-v|| \leq \delta \mbox{ \ , for some \ } v\in S \} ,
$$
for $ c, \rho, \delta \in (0, + \infty ) $. Since the proofs of Lemma \ref{l118} and Corollary \ref{l119} can be done as \cite[Lemma 4.6]{[6]}, we omit it here.

\bl\label{l118}
Let $ c, \rho > 0 $. Then, there exists $ \delta_0 = \delta(c, \rho)>0 $ such that, if $ \delta \in (0, \delta_0) $, then $ ||I'(u)||_{X'}(1+||u||_X) \geq 8 \delta $, for all $ u\in \tilde{S}_{2\delta} $ with $ I(u)\in [c-2\delta^2 , c+ 2 \delta^2] $.
\el

\bc\label{l119}
Let $ c, \rho > 0 $. Then, there exists $ \tilde{\delta}_0 = \tilde{\delta}_0(c, \rho)>0 $ such that, if $ \delta \in (0, \tilde{\delta}_0) $, then $ ||I'(u)||_{X'}(1+||u||_X) \geq 8 \delta $, for all $ u\in S_{2\delta} $ with $ I(u)\in [c-2\delta^2 , c+ 2 \delta^2] $.
\ec

\bl\label{l120}
Let $ c> 0 $. Then, there exists $ \rho_1 = \rho_1 (c) > 0 $ such that, for all $ \rho \in (0, \rho_1)  $, we have

\noindent \textbf{(i)} $ A_{c, \rho} \cap D = \emptyset $ .

\noindent \textbf{(ii)} There exists $ \varepsilon = \varepsilon(c, \rho) >0 $ and a function $ \phi : I^{c+ \varepsilon}\setminus A_{c, \rho} \RA I^{c-\varepsilon} $, continuous and even, such that $ D \subset I^{c-\varepsilon} $ and $ \phi\big\vert_{D} = id $.
\el
\begin{proof}
Item (i) follows from the application of Proposition \ref{p3} and Lemma \ref{l12}, arguing by contraction. For item (ii), let $ \delta_0 $ as given by Lemma \ref{l119} and $ \delta\in (0, \delta_0) $, such that $ \delta^2 < \frac{c}{2} $. Take $ \varepsilon = \delta^2 $. Then, from Deformation's Lemma 2.6 of \cite{[def]}, there exists $ \eta : [0, 1]\times X \RA X $, continuous, satisfying

\noindent \textbf{(a)} $ \eta(t, u)=u $, if $ t=0 $ or $ u\not\in I^{-1}([c-2\varepsilon , c+ 2 \varepsilon])\cap S_{2 \delta} $ ;

\noindent \textbf{(b)} $ \eta(1, I^{c+ \varepsilon}\cap (X \setminus A_{c, \rho})) \subset I^{c-\varepsilon} $;

\noindent \textbf{(c)} $ t \mapsto I(\eta(t, u)) $ is non-increasing, for all $ u\in X $.

Moreover, since $ I $ is odd, it is possible to modify the proof of the refereed lemma, such as in \cite{[defR]}, to obtain as well

\noindent \textbf{(d)} $ \eta(t, -u) = - \eta(t, u) $, for all $ t \in [0, 1] $ and $ u\in X $.

Define $ \phi : I^{c+ \varepsilon}\setminus A_{c, \rho} \RA I^{c-\varepsilon}  $ by $ \phi(u)=\eta(1, u) $. Note that item (b) is equivalent to $ \eta(1, I^{c+ \varepsilon}\cap \setminus A_{c, \rho}) \subset I^{c-\varepsilon} $, which guarantee that $ \phi $ is well-defined. Also, as $ \eta $ is continuous, $ \phi $ is continuous as well, and from item (d), $ \phi $ is even. 

Moreover, once $ \varepsilon = \delta^2 < \frac{c}{2} $, $ 0\not\in [ c-2 \varepsilon , c+ 2 \varepsilon] $. Thus, $ D \cap I^{-1}([ c-2 \varepsilon , c+ 2 \varepsilon]) = \emptyset $. Hence, from item (a), if $ u\in D $, $ \phi(u)=\eta(1, u)=u $, which implies that $ \phi\big\vert_{D}=id $. 

Finally, if $ u\in D $, $ I(u)\leq 0< c-\varepsilon $, then $ u\in I^{c-\varepsilon} $. 
\end{proof}

\bp\label{p112}
Let $ k\in \N $. Then, $ c_k $ is a critical value of $ I $.
\ep
\begin{proof}
Arguing by contradiction, from Lemma \ref{l120} and from items (ii) and (iii) of Lemma \ref{l19}, one get that $ c_k $ is a critical value of $ I $.
\end{proof}

\bp\label{p113}
We have $ c_k \RA + \infty $, as $ k \RA + \infty $.
\ep
\begin{proof}
Suppose that there exists $ M >0 $, such that $ c_k < M $, for all $ k\in \N $. From Remark \ref{obs15}, $c_k$ is monotonically nondecreasing. Then, there exists $ c > 0 $ such that $ c_k \RA c $.

From Proposition \ref{p19}, there exists $ \rho_0 > 0 $ such that $ \gamma(A_{c, \rho}) < + \infty $, for all $ \rho \in (0, \rho_0) $. Also, from items (i) and (ii) of Lemma \ref{l120} and items (ii) and (iv) of Lemma \ref{l19}, there exists $ \varepsilon >0 $ such that
$$
\gamma_D (I^{c+ \varepsilon}) = \gamma_D((I^{c+ \varepsilon}\setminus A_{c, \rho})\cup A_{c, \rho}) \leq \gamma_D (I^{c+ \varepsilon}\setminus A_{c, \rho}) + \gamma(A_{c, \rho}) < + \infty .
$$
Once more, as $ c_k \RA c $ monotonously nondecreasing and by definition of values $ c_k $, we have $ \gamma_D (I^{c+ \varepsilon}) \RA + \infty $, as $ k \RA + \infty $, which leads to a contradiction. 

Therefore, $ c_k \RA + \infty $, as $ k \RA + \infty $.
\end{proof}

\begin{proof}[Proof of Theorem \ref{t2}] 
From Proposition \ref{p113}, we can extract a subsequence of $ (c_k) $ such that $ c_k \RA + \infty $ monotonously increasing. Then, from Proposition \ref{p112}, there exists $ \uk \in X $ satisfying $ I(\uk)=c_k $ and $ I'(\uk)=0 $, for all $ k\in \N $. Since $ (c_k) $ is monotone, $ c_i \neq c_j $, when $ i\neq j $, and $ c_k > c_1 >0 $, for all $ k\geq 2 $. Then, we conclude that the functions $ u_k $ are distinct and that $ \uk \neq 0 $, for all $ n\in \N $. Also, as $ I $ is odd, the same holds for $ -\uk $ and $ I(\pm \uk) \RA + \infty $. 
\end{proof}

\noindent \textbf{Data availability statement:} The data that supports the findings of this study are available within the article [and its supplementary material].
\\

\noindent \textbf{Acknowledgements:} This work was perfomed and completed during the  first author PhD graduate course at Federal University of São Carlos.

\end{document}